\documentclass{amsart}
\usepackage{geometry} 
\usepackage{amsfonts}
\usepackage{amsmath}
\usepackage{graphicx}
\usepackage{mathrsfs}
\usepackage{amssymb}
\usepackage{amsthm}
\usepackage{tikz}
\usepackage{tikz-cd}
\usepackage{listings}
\usepackage{color}
\usepackage{enumitem}
\usepackage{imakeidx}
\usepackage[colorinlistoftodos]{todonotes}
\usepackage{upgreek}
\usepackage[pagebackref=true]{hyperref} 
\usepackage{mdframed}
\usepackage{pgfplots}
\usepackage{manfnt}
\usepackage[new]{old-arrows}
\usepackage{changepage}
\usepackage{stmaryrd}
\usepackage{appendix}
\usepackage[normalem]{ulem}
\usepackage{ytableau}
\usepackage{blkarray}
\usepackage{mathtools}
\usepackage{verbatim}

\DeclareMathAlphabet{\matheuler}{U}{eus}{m}{n}

\pgfplotsset{width=7cm,compat=1.9}

\makeindex[intoc]

\definecolor{codegreen}{rgb}{0,0.6,0}
\definecolor{codegray}{rgb}{0.5,0.5,0.5}
\definecolor{codepurple}{rgb}{0.58,0,0.82}
\definecolor{backcolour}{rgb}{0.95,0.95,0.92}

\lstdefinestyle{mystyle}{
	backgroundcolor=\color{backcolour},   
	commentstyle=\color{codegreen},
	keywordstyle=\color{magenta},
	numberstyle=\tiny\color{codegray},
	stringstyle=\color{codepurple},
	basicstyle=\footnotesize,
	breakatwhitespace=false,         
	breaklines=true,                 
	captionpos=b,                    
	keepspaces=true,                 
	numbers=left,                    
	numbersep=5pt,                  
	showspaces=false,                
	showstringspaces=false,
	showtabs=false,                  
	tabsize=2
}

\allowdisplaybreaks

\newcommand{\MAT}[7]{\begin{pmatrix}{#1}&{#2}&{#3}&{#4}\\{#5}&{#6}&{#7}&\MATHelper}
\newcommand{\MATHelper}[9]{{#1}\\{#2}&{#3}&{#4}&{#5}\\{#6}&{#7}&{#8}&{#9}\end{pmatrix}}

\DeclareMathOperator{\Hom}{Hom}

\DeclareMathOperator{\Res}{Res}

\DeclareMathOperator{\Nm}{Nm}

\DeclareFontFamily{U}{wncy}{}
\DeclareFontShape{U}{wncy}{m}{n}{<->wncyr10}{}
\DeclareSymbolFont{mcy}{U}{wncy}{m}{n}
\DeclareMathSymbol{\Sha}{\mathord}{mcy}{"58} 

\newcommand{\val}{\operatorname{val}}

\makeatletter
\newcommand{\bigperp}{%
  \mathop{\mathpalette\bigp@rp\relax}%
  \displaylimits
}

\newcommand{\bigp@rp}[2]{%
  \vcenter{
    \m@th\hbox{\scalebox{\ifx#1\displaystyle2.1\else1.5\fi}{$#1\perp$}}
  }%
}
\makeatother

\DeclareMathOperator{\Lie}{Lie}

\DeclareMathOperator{\Spec}{Spec}

\let\templim\lim
\renewcommand{\lim}{\templim\limits}

\DeclareMathOperator{\cyc}{cyc}
\DeclareMathOperator{\anti}{anti}

\newcommand\homses[4]{
    \begin{tikzcd}[ampersand replacement=\&]
        0\ar[r]\&{\displaystyle #1}\ar[r, "{#2}"]\ar[r]\&{\displaystyle #3}\ar[r, "{#4}"]\&
        \homsesHELPER
}
\newcommand\homsesHELPER[9]{
        {\displaystyle #1}\ar[r]\& 0\\
        0\ar[r]\&{\displaystyle #5}\ar[r, "{#6}"]\&{\displaystyle #7}\ar[r, "{#8}"]\&{\displaystyle #9}\ar[r]\& 0
        \ar[from=lllu, to=lll, "{#2}"]\ar[from=llu, to=ll, "{#3}"]\ar[from=lu, to=l, "{#4}"]
    \end{tikzcd}
}

\renewcommand{\phi}{\varphi}
\newcommand{\mfp}{\mathfrak p}

\newcommand{\mf}{\mathfrak}

\newcommand{\mfq}{\mf q}

\newcommand{\mc}{\mathcal} 

\newcommand{\mcO}{\mc O}

\newcommand{\mcR}{\mc R}

\newcommand{\mcU}{\mc U}

\newcommand{\mcE}{\mc E}

\newcommand{\mcS}{\mc S}

\newcommand{\mcC}{\mc C}

\newcommand{\F}{\mathbb F}
\newcommand{\Q}{\mathbb Q}
\newcommand{\Z}{\mathbb Z}

\newcommand{\A}{\mathbb A}
\newcommand{\G}{\mathbb G}
\newcommand{\K}{\mathbb K}

\renewcommand{\tau}{\uptau}
\renewcommand{\P}{\mathbb P}

\newcommand{\dparens}[1]{\!\left(\!\left(#1\right)\!\right)}
\newcommand\dps\dparens 

\newcommand{\into}{\hookrightarrow}

\makeatletter
\newcommand*{\da@rightarrow}{\mathchar"0\hexnumber@\symAMSa 4B }
\newcommand*{\da@leftarrow}{\mathchar"0\hexnumber@\symAMSa 4C }
\newcommand*{\xdashrightarrow}[2][]{%
  \mathrel{%
    \mathpalette{\da@xarrow{#1}{#2}{}\da@rightarrow{\,}{}}{}%
  }%
}
\newcommand{\xdashleftarrow}[2][]{%
  \mathrel{%
    \mathpalette{\da@xarrow{#1}{#2}\da@leftarrow{}{}{\,}}{}%
  }%
}
\newcommand*{\da@xarrow}[7]{%
  \sbox0{$\ifx#7\scriptstyle\scriptscriptstyle\else\scriptstyle\fi#5#1#6\m@th$}%
  \sbox2{$\ifx#7\scriptstyle\scriptscriptstyle\else\scriptstyle\fi#5#2#6\m@th$}%
  \sbox4{$#7\dabar@\m@th$}%
  \dimen@=\wd0 %
  \ifdim\wd2 >\dimen@
    \dimen@=\wd2 %
  \fi
  \count@=2 %
  \def\da@bars{\dabar@\dabar@}%
  \@whiledim\count@\wd4<\dimen@\do{%
    \advance\count@\@ne
    \expandafter\def\expandafter\da@bars\expandafter{%
      \da@bars
      \dabar@ 
    }%
  }%
  \mathrel{#3}%
  \mathrel{%
    \mathop{\da@bars}\limits
    \ifx\\#1\\%
    \else
      _{\copy0}%
    \fi
    \ifx\\#2\\%
    \else
      ^{\copy2}%
    \fi
  }%
  \mathrel{#4}%
}
\makeatother
\newcommand\xdashto\xdashrightarrow
\newcommand\xdashfrom\xdashleftarrow

\makeatletter
\providecommand{\leftsquigarrow}{%
  \mathrel{\mathpalette\reflect@squig\relax}%
}
\newcommand{\reflect@squig}[2]{%
  \reflectbox{$\m@th#1\rightsquigarrow$}%
}
\makeatother

\newcommand*{\rom}[1]{\textup{\uppercase\expandafter{\romannumeral#1}}}

\renewcommand{\bar}{\overline}

\renewcommand\t\text 
\newcommand\ttt\texttt

\newcommand{\important}[1]{\textit{#1}}
\newcommand\noteworthy\important

\pgfmathdeclarerandomlist{thuses}{{Thus}{Therefore}{As such}{As we wished}{By this great fortune}{Alors}{Donc}{Following this thread}{Consequently}{Thusforth}{Rejoice}{This can only mean one thing...}{Ergo}{Voila}}

\pgfmathdeclarerandomlist{fibers}{{fiber}{fibre}}

\pgfmathdeclarerandomlist{wewins}{{we win}{it's over}{all is good}{victory has been achieved}{we've crossed the finish line}{we're happy campers}{we've done it}{that's a wrap}{we have the high ground}}

\tikzset{%
	symbol/.style={%
		,draw=none
		,every to/.append style={%
			edge node={node [sloped, allow upside down, auto=false]{$#1$}}}
	}
}

\makeatletter
\def\renewtheorem#1{%
	\expandafter\let\csname#1\endcsname\relax
	\expandafter\let\csname c@#1\endcsname\relax
	\gdef\renewtheorem@envname{#1}
	\renewtheorem@secpar
}
\def\renewtheorem@secpar{\@ifnextchar[{\renewtheorem@numberedlike}{\renewtheorem@nonumberedlike}}
\def\renewtheorem@numberedlike[#1]#2{\newtheorem{\renewtheorem@envname}[#1]{#2}}
\def\renewtheorem@nonumberedlike#1{  
	\def\renewtheorem@caption{#1}
	\edef\renewtheorem@nowithin{\noexpand\newtheorem{\renewtheorem@envname}{\renewtheorem@caption}}
	\renewtheorem@thirdpar
}
\def\renewtheorem@thirdpar{\@ifnextchar[{\renewtheorem@within}{\renewtheorem@nowithin}}
\def\renewtheorem@within[#1]{\renewtheorem@nowithin[#1]}
\makeatother

\DeclareDocumentEnvironment{MyFrame}{O{1cm}O{0.4pt}O{0.8cm}O{black}O{3}O{2ex}}
{\par\hfill\rlap{%
		\bgroup\color{#4}%
		\hskip-\dimexpr#1-#3\relax\rule{#1}{#2}%
		\hskip-\dimexpr#1/#5\relax\rule[-\dimexpr#1-\dimexpr#1/#5\relax]{#2}{#1}%
		\egroup
	}%
	\vskip-\dimexpr#1/#5+\dimexpr#1/#5-#6\relax%
}
{\par\nobreak\offinterlineskip\vskip-\dimexpr#1/#5+\dimexpr#1/#5-#6\relax\noindent%
	\hskip-#3\bgroup\color{#4}%
	\rule{#1}{#2}\hskip-\dimexpr#1-\dimexpr#1/#5-#2\relax%
	\rule[-\dimexpr#1/#5-#2\relax]{#2}{#1}\egroup\par
}

\newtheorem{innercustomgeneric}{\customgenericname}
\providecommand{\customgenericname}{}
\newcommand{\newcustomtheorem}[2]{%
	\newenvironment{#1}[1]
	{%
		\renewcommand\customgenericname{#2}%
		\renewcommand\theinnercustomgeneric{##1}%
		\innercustomgeneric
	}
	{\endinnercustomgeneric}
}

\theoremstyle{break}
\newtheorem{algorithm}{Algorithm}
\newtheorem{thm}{Theorem}[section]

\newtheorem{lemma}[thm]{Lemma}

\renewtheorem{conj}[thm]{Conjecture}
\newtheorem{propn}[thm]{Proposition}

\newcustomtheorem{Prob}{Problem}
\newcustomtheorem{Exc}{Theoretical Exercise}
\newtheorem{example}[thm]{Example}

\theoremstyle{definition}
\newtheorem{defninner}[thm]{Definition}
\newtheorem{cond}[thm]{Condition}

\newtheorem{recinner}[thm]{Recall}

\newtheorem*{simpassump}{Simplifying Assumption}

\newtheorem*{nonexinner}{Non-example}
\newtheorem{warninner}[thm]{Warning}

\newtheorem*{ansinner}{Answer}

\newtheorem*{histinner}{History}

\theoremstyle{remark}
\newtheorem{rmk}[thm]{Remark}

\newcommand\remsymbol{$\circ$}
\newcommand\anssymbol{$\star$}
\newcommand\warnsymbol{$\bullet$}
\newcommand\proofsymbol{$\blacksquare$}
\newcommand\nonexsymbol{$\triangledown$}
\newcommand\recsymbol{$\odot$}
\newcommand\defnsymbol{$\diamond$}
\newcommand\histsymbol{$\ominus$}

\newenvironment{defn}[1][]{\begin{defninner}[#1]\pushQED{\qed}\renewcommand\qedsymbol\defnsymbol}{\popQED\end{defninner}}

\let\proofinner=\proof
\newcommand\myproof[1][Proof]{\proofinner[#1]\renewcommand\qedsymbol\proofsymbol}
\def\proof{\myproof}
\DeclareMathOperator{\red}{red}

\begin{document}

\newcommand{\aash}[1]{{\color{blue} \sf $\clubsuit\clubsuit\clubsuit$ Aash: [#1]}}

\title{Finding Integral Points on Elliptic Curves over Imaginary Quadratic Fields}
\author[Aashraya Jha]{Aashraya Jha}
\address{aashjha@bu.edu\\ Department of Mathematics and Statistics, 665 Commonwealth Avenue, Boston MA 02215 }

\begin{abstract}
  We determine the quadratic Chabauty set for integral points on elliptic curves of rank $2$ defined over imaginary quadratic fields using quadratic Chabauty. This builds on the work of Bianchi \cite{bianchi2020integral} and Balakrishnan et al. \cite{balakrishnan2021explicit}. We give the first instance of the implementation of anticyclotomic heights for curves which are not base changes, along with an implementation of a certain sieve for elliptic curves introduced by Balakrishnan et al. \cite{balakrishnan2017computing} and used by Bianchi \cite{bianchi2020integral} to determine integral points of rank $2$. We give the first example of the determination of the integral points of an elliptic curve of rank $2$ defined over an imaginary quadratic field, which is not a base change via quadratic Chabauty.
\end{abstract}

\date{November 2, 2023}

\maketitle
\tableofcontents

\section{Introduction}\label{S:intro}

Let $K$ be a number field with $O_K$ its ring of integers. Let $C/K$ be a smooth, projective, and geometrically irreducible curve, henceforth called a \emph{nice curve}. Let $\mcU/O_K$ be an affine scheme with projective closure $\mcC$  such that $\mcC\times \Spec K\cong C$. Let $g$ be the genus of $C$. In $1929$, Siegel \cite{siegel2014einige} showed that the set $\mcU(O_K)$ is finite if $g\geq 1$. For $g\geq 2$, this was superseded by Faltings's theorem in 1983 \cite{faltings1983endlichkeitssatze}, which shows $C(K)$ is finite for $g\geq 2$.

Neither of these proofs can be made effective; that is, we can not use them to explicitly determine the sets $C(K)$ or $\mcU(O_K)$. There has been significant progress in determining the sets $\mcU(O_K)$ and $C(K)$ via $p$-adic methods of Chabauty and Coleman. Let $J$ be the Jacobian of $C$ and $r$ be the rank of  $J(K)$.
 In his paper \cite{chabauty1941points}, Chabauty showed the set $C(K)$ is finite if $r<g$.

We give a heuristic explaining why we might expect this when $K=\Q$. We first pick a prime of good reduction, say $p$.

We can identify $C$ as a subvariety of $J$ via the Abel-Jacobi map $\iota:C\into J$. If we have $r\leq g-1$, then the closure $\overline{J(K)}\subseteq J(\Q_p)$ has dimension at most $r$ (as an analytic variety over $\Q_p$), and the curve $C_{\Q_p}$ has dimension $1$. Since $r+1\leq g$, we expect the intersection $\overline{J(K)}\cap C_{\Q_p}$ to be finite. 

Chabauty constructs functions that vanish on this intersection, which Coleman \cite{coleman1985torsion} later identified as $p$-adic (Coleman) integrals of holomorphic differential forms. The method of Chabauty-Coleman relies on the image of the map
\begin{equation}\label{eqn: Chabauty-Coleman}
    \log: J(\Q) \otimes \Q_{p} \to H^0(X_{\Q_p}, \Omega^1)^{\vee}
\end{equation}
having positive codimension, which is used to write down abelian integrals vanishing
at the rational points.

Kim's work \cite{kim2009unipotent} gave rise to the \emph{non-abelian} Chabauty program, which aims to remove the restriction of $r<g.$ Kim's method considers unipotent quotients of the $\Q_p$-\'etale fundamental group of $C$ and looks at the Selmer schemes attached to these quotients. In \cite{kim2010massey} and its appendix \cite{balakrishnan2011appendix}, the authors show that a certain $p$-adic locally analytic function is constant on integral points of elliptic curves with Tamagawa product $1$. This $p$-adic function is a linear combination of an iterated Coleman integral and the square of an abelian Coleman integral.

Let $p$ be a prime of good reduction. We let \[h\colon C(\Q)\to \Q_p\] denote the global height of Coleman and Gross \cite{coleman1989p}, and for finite places $v$ of $K$, we let \[h_v\colon C(\Q_v)\to\Q_p\] denote the local heights of Coleman and Gross. We have $h=\sum_v h_v|_{C(\Q)}$. In \cite{balakrishnan2015coleman}, the authors extend the function $h-h_p$ to a locally analytic function on $C(\Q_p)$, and show that this function is a scalar multiple of the Coleman function in \cite{kim2010massey} and \cite{balakrishnan2011appendix} . In \cite{balakrishnan2016quadratic}, the authors further extend this approach to hyperelliptic curves. Let $C$ be a hyperelliptic curve defined by $y^2-f(x)$ where $f$ has degree $2g+1$, and let \[\mcU\coloneqq\Spec\Z[x,y]/(y^2-f(x))).\]
They prove the following theorem regarding integral points on hyperelliptic curves.
\begin{thm}\cite[Theorem 3.1]{balakrishnan2016quadratic}
    Suppose that $r = g$ and that the log map in
\eqref{eqn: Chabauty-Coleman} is an isomorphism. Then there exists an explicitly computable finite set $T \subset \Q_p$ and an explicitly computable non-constant Coleman function $\rho:\mcU(\Z_p)\to \Q_p$ such that $\rho(\mcU(\Z)) \subseteq T$. 
\end{thm}

We would like to extend this method to curves over more general number fields $K$. Based on an idea of Wetherell, Siksek \cite{siksek2013explicit} looks at the restriction of scalars of $V\coloneqq \Res^K_{\Q}C$ and $A\coloneqq \Res^K_{\Q}J$. He exploits the fact $C(K)=V(\Q)$ and that
\begin{equation}\label{eqn:Weil restriction}
    V(\Q)\subseteq V(\Q_p)\cap \bar{A(\Q)}\subseteq A(\Q_p).
\end{equation}

Let $B\coloneqq V(\Q_p)\cap \bar{A(\Q)}$. If $r\leq d(g-1)$, one expects $B$ to be finite due to reasons of dimension. 

\begin{rmk}
Siksek \cite{siksek2013explicit} provides a statement that he thinks is \textbf{possibly correct}. Let $L\subseteq K$ be a subfield and let $D/L$ be a nice curve such there is an isomorphism $D\times_{\Spec L} \Spec K\cong C$. Let $J^D$ be the Jacobian of $D$ and let $r_{D}$ be the rank of $J^D(L)$. Then if for all $L\subseteq K$ and all $D/L$ if
\begin{equation}\label{eqn: Chab over NF}
    r_D\leq [L:\Q](g-1),
\end{equation}
then $B$ is finite. Dogra \cite{dogra2023unlikely} shows this statement is not true; he constructs a hyperelliptic genus $3$ curve for which the the set $B$ is infinite even though it satisfies \eqref{eqn: Chab over NF} for each such $D/L$. On the other hand, he shows if \eqref{eqn: Chab over NF} holds and further if $J$ does not share a component with any of its conjugates over $\bar{\Q}$, we have finiteness of the set $B$.
\end{rmk}

In the paper \cite{balakrishnan2021explicit}, the authors exploit the idea of using Weil restrictions and multiple $p$-adic height functions to give a quadratic Chabauty method over number fields. They consider the map 
\begin{equation}\label{eqn: Quadratic Chabuty}
    \log: J(K)\otimes \Q_p \to (\Res^K_{\Q}J)(\Q)\otimes \Q_p \to \Lie(\Res^K_{\Q}(J))_{\Q_p}.
\end{equation}

Let $p$ be a prime which splits completely in $K$. For $1\leq i\leq m$, let $\psi_i:K\into \Q_p$ denote all the embeddings of $K$ into $\Q_p$.
Let \[O_K\otimes \Z_p = \prod_{i=1}^m O_K\otimes_{\psi_i}\Z_p \text{ and } \mcU(O_K\otimes \Z_p)\coloneqq \prod_{i=1}^m \mcU(O_K\otimes_{\psi_i}\Z_p).\]

Also let $\psi:\mcU(O_K)\into \mcU(O_K\otimes \Z_p)$ denote the natural embedding induced by the $\psi_i$ for $1\leq i
\leq m$. We refer the reader to Section \ref{Section:idele class characters} for the definition of $p$-adic id\`ele class characters. They prove the following theorem:

\begin{thm}[\cite{balakrishnan2021explicit}] \label{thm:Intergal points number fields}
Let $p$ be a prime such that $\mcU$ has good reduction at all primes above $p$, and let $\chi$ be a $p$-adic id\`ele class character.
 Suppose the map in \eqref{eqn: Quadratic Chabuty} is injective. Then there exists an explicitly computable finite set $T^{\chi}\subset \Q_p$  and an explicitly computable non-constant locally analytic
function $\rho^{\chi} : \mcU(O_K\otimes \Z_p) \to \Q_p$, both dependent on $\chi$, such that $\rho^{\chi}(\psi(\mcU(O_K)))\subseteq T^{\chi}$.
    
\end{thm}

For practical applications, computing the set of integral points via this method has been limited to quadratic number fields for curves defined over $\Q$. Some of the difficult steps in extending this method to higher degree number fields include finding the generators of (a finite-index subgroup of) the Mordell--Weil group $J(K)$ when the curve is defined over a larger degree number field, implementing Coleman integration and computing heights over field extensions of $\Q_p$.

We review examples that have been computed in the literature thus far. In \cite{balakrishnan2021explicit}, the authors compute the set of $\Z[\zeta_3]$-points of a genus $2$ curve defined over $\Q$. The \emph{quadratic Chabauty set} is the set of is the preimage of $T^{\chi}$ under $\rho^{\chi}$ as $\chi$ ranges over linearly independent id\`ele class characters. In particular, it is a set of $p$-adic points in $\mcU(O_K\otimes\Z_p)$. The authors of \cite{balakrishnan2021explicit} compute finite quadratic Chabauty sets, which are supersets of the $O_K$-points, of two elliptic curves where $K$ is $\Q(\sqrt{5})$ and $\Q(\sqrt{3})$. To find these sets, they use the \emph{cyclotomic height} and a relation among the linear forms given by abelian Coleman integrals.

We note that the Mordell--Weil sieve cannot be applied to get rid of auxiliary $p$-adic points which are zeroes of Coleman functions on elliptic curves, since the sieve crucially relies on the fact that the curve embeds as proper subvariety of its Jacobian and consequently the map 
\[C_{\F_q}(\F_q)\to J_{\F_q}(\F_q)\]
 is not surjective. 

Gajovi\'c and M\"uller \cite{gajovic2023linear} compute the $\Z[\sqrt{7}]$-points of a hyperelliptic curve which is not a base change using \emph{linear quadratic Chabauty}. Their method applies to a hyperelliptic curve in the form \[y^2=f(x),\]
\noindent where the degree of $f$ is even and $f$ is monic. This complements the work in 
 \cite{balakrishnan2021explicit}, where the authors provide a method when the degree is odd. We remark the linear quadratic Chabauty method of Gajovi\'c and M\"uller uses an innovation in \cite{gajovic2023computing} to compute $p$-adic heights of hyperelliptic curves. Based on the work of \cite{balakrishnan2012computing}, they provide an algorithm which cleverly uses the degree $0$ divisor above infinity to compute heights on hyperelliptic curves, which is significantly faster if $\deg(f)$ is odd.
 
    Bianchi \cite{bianchi2020integral} looks at the Mordell curve with \href{https://beta.lmfdb.org/EllipticCurve/2.0.3.1/20736.1/CMd/1}{LMFDB label 20736.1-CMd1} and minimal Weierstrass equation
\[y^2=x^3-4\]
base changed to $\Q(\zeta_3)$,
and uses the theory of Mazur-Stein-Tate \cite{mazur2006computation} heights on elliptic curves. This heavily relies on the $p$-adic sigma function developed by Mazur and Tate \cite{mazur1991p}. 
On computeing the quadratic Chabauty locus, she then employs a sieve that exploits the size of the set of mod $7$ and mod $13$ points of $E$. This allowed her to determine all the $\Z[\zeta_3]$ points of this Mordell curve.

In this paper, we give a method to compute the integral points of rank $2$ elliptic curves over imaginary quadratic fields which are \emph{not} base changes. We follow the strategy proposed by Bianchi in \cite{bianchi2020integral}. Some of the simplifications both in the quadratic Chabauty method and the application of the sieve are lost when one works with a curve that is not a base change. Our work gives the first instance of a calculation of the \emph{anticyclotomic} height of an elliptic curve, which is not a base change from $\Q$, see Proposition \ref{Propn: Anticyclotomic height} and Algorithm \ref{algo:heights}. We also present the first example where the integral points of an elliptic curve over an imaginary quadratic field which is not a base change of an elliptic curve over $\Z$ have been determined via the method of quadratic Chabauty. 

The most common algorithms for computing the set of integral points on elliptic curves are based on the study of elliptic logarithms. The implementation in Magma is based on an algorithm of Stroeker and Tzanakis \cite{stroeker1994solving} which uses linear (group) relations between integral points and the generators of the free component of the Mordell--Weil group transformed into a linear form in elliptic logarithms. This was extended to number fields in the work of Smart and Stephens \cite{smart1997integral}. This implementation seems to be only available for elliptic curves defined over totally real fields. Sage currently has no implementation for computing the set of integral points over number fields apart from the rationals.  We thus give a method to determine integral points for elliptic curves over imaginary quadratic fields.

We can use our quadratic Chabauty method to prove the following result:
\begin{thm}
Let $O_K=\Z[\zeta_3]$ with fraction field $K$. 
Consider the scheme $\mcU_1/O_K$ cut out by the Weierstrass equation 
    \[y^2+(\zeta_3+2)y=x^3+(-\zeta_3-2)x^2+(\zeta_3+1)x.\]

This curve has \href{https://beta.lmfdb.org/EllipticCurve/2.0.3.1/134689.3/CMa/1}{LMFDB label 134689.3-CMa1} Then \begin{align*}
    \mcU_1(O_K)=&\{(-3 : -8\zeta_3 - 4), (-3 : 7\zeta_3 +2 ),(4\zeta_3+4 : -8\zeta_3 - 4), (0 : 0 ), (0 : -\zeta_3 - 2), (\zeta_3+1 : 0 ), \\(\zeta_3+1 : -\zeta_3 - 2)
    & (4\zeta_3+4 : 7\zeta_3 +2 ), (1 : 0), (1 : -\zeta_3 - 2), (-3\zeta_3 + 1 : -8\zeta_3 - 4), (-3\zeta_3 + 1 : 7\zeta_3 + 2)\}.
    \end{align*} 
\end{thm}

\begin{rmk}
    The data in the LMFDB is presented in terms of the algebraic integers $a=\zeta_6$. Hence the equations on LMFDB are obtained by replacing $\zeta_3$ above with $a-1$.
\end{rmk}

The paper is organised as follows. We look at the theory of heights as formulated by Mazur, Stein, and Tate \cite{mazur2006computation} in Section \ref{sec: Heights}. We outline the strategy to determine the quadratic Chabauty sets for a given elliptic curve in Section \ref{Sec:QC strategy}. We then look at a sieve to eliminate spurious $p$-adic points in Section \ref{sec: sieve}.
Finally, we compute explicit examples in Section \ref{Sec:Examples}. 

While we can determine the quadratic Chabauty set for elliptic curves over imaginary quadratic fields of class number $1$, this method is not very effective in determining the set of integral points. The major stumbling block is the sieve we use. 
To use the sieve, we need the curve to have trivial torsion, and we need the cardinalities of reductions at two split primes to satisfy a particular condition; see Condition \ref{condition: reduction}. Unfortunately, this condition is not satisfied too often, and when it is satisfied, the sieve often does not end up eliminating all mock integral points. For example, we found $3528$ curves over $\Q(\zeta_3)$ of rank $2$, trivial torsion, which were not base changes. Of these, only $36$ satisfied Condition \ref{condition: reduction}. Amongst these curves, the sieve described in Section \ref{sec: sieve} only managed to eliminate all mock integral points in one case. Thus, developing a better sieving algorithm for $p$-adic points of elliptic curves would help improve this method significantly. 

The code used in this paper is available on \href{https://github.com/AashrayaJha/QC_ECIm}{Github} \cite{gitrepo}. 

\section*{Acknowledgements}

We would like to thank Jennifer Balakrishnan for her constant support during this project and all the helpful comments on the write-up. We would also like to thank Alex Best, Alex Betts, Netan Dogra and Steffen M\"uller for illuminating conversations. We are grateful for support from the Simons Foundation as part of the Simons Collaboration on Arithmetic Geometry, Number
Theory, and Computation \#550023
and NSF grant DMS-1945452. 
\section*{Notation}\label{S:notn}
\begin{alignat*}{2}
&K  &&\text{A number field},\\
&O_K     &&\text{The ring of integers of $K$},\\
&\mfp,\mfq,v    &&\text{Finite places of $K$}, \\    
&K_v    &&\text{The completion of $K$ along a place $v$},\\
&\psi_i:K\into \Q_p &&\text{Embeddings of $K$ into $\Q_p$ for completely split primes $p$.}\\
& O_K\otimes \Z_p\coloneqq \prod_{i=1}^m O_K\otimes_{\psi_i}\Z_p &&\text{The product of tensor products given by all embeddings. } \\
& A_K^{\times}&& \text{The group of id\`eles of a number field.} \\
&E/K    &&\text{An elliptic curve over $K$},\\
&E_{K_v} \coloneqq E\times_{\Spec K}\Spec K_v &&\text{ The base change of $E$ to $K_v,$}\\
&\mcE/O_K &&\text{An integral model of $E$, often the minimal Weierstrass model,}\\
&\mcU/O_K\coloneqq \mcE\setminus \{\mcO\} &&\text{$\mcE$ minus the identity section}\\
&E_{\F_{\mfp}}\coloneqq \mcE\times O_K/\mfp && \text{The reduction of $E$ at $\mfp$}. \\ 
\end{alignat*}

\begin{rmk}
For a rational prime $p$  which splits completely in $K$ and for a prime $\mfp_i\mid p$ of $K$, we will also use $\psi_i$ to denote the isomorphism \[\psi_i\colon K_{\mfp_i}\to\Q_p\]
which extends the embedding $\psi_i\colon K\into \Q_p$.
\end{rmk}

\begin{rmk}\label{Defn: diagonal embedding}
    We set \[ 
        \psi\colon K \to\prod_{\mfp|p} K_{\mfp}, \,
        x \mapsto \prod \psi_i(x),
\]
to be the diagonal embedding of $K$ into its completions at places above $p$.
\end{rmk}
\section{Computing non-Archimedean heights on elliptic curves} \label{sec: Heights}

A fundamental tool used in quadratic Chabauty is the theory of $p$-adic heights. In this section, after recalling some basic properties of heights, we will focus on the computational aspects of $p$-adic heights of elliptic curves over imaginary quadratic fields.

We shall use the definition of $p$-adic heights from \cite{mazur1983canonical} and \cite{mazur2006computation}. 

In \cite{mazur1983canonical} the authors define height parings valued in various abelian groups for abelian varieties defined over various fields, but we stick to the case where $E$ is an elliptic curve over a number field $K$, and our height pairing is valued in $\Q_p$, for some rational prime $p$. Given an (admissible) id\`ele  class character \footnote{See Definition \ref{def: Idele class character}} \[\chi:\A_K^{\times}/K^{\times}\to\Q_p\] they show there exists a ``canonical''  pairing 
\[(\cdot,\cdot)_{\chi}:E(K)\times E(K)\to \Q_p\] which is symmetric and bilinear and therefore is a height pairing. This height satisfies some nice functorial properties, see \cite[Section 4]{mazur1983canonical}.

In \cite{mazur2006computation}, the authors give explicit formulas for the canonical pairings described in \cite{mazur1983canonical} as a sum of local pairings and an algorithm to compute these pairings. David Harvey \cite{harvey2008efficient} provides a more efficient algorithm to calculate the height associated to the cyclotomic character for elliptic curves defined over $\Q$, which we modify to obtain cyclotomic heights in the number field case.

The main inputs in computing $p$-adic heights on elliptic curves are id\`ele class characters, the  $p$-adic sigma function introduced in \cite{mazur1991p}, and denominators of the $x$-coordinate of given points. We say more about these in the upcoming subsections.

 \begin{simpassump}
     We assume henceforth $K$ has class number $1$ for ease of exposition and to speed up computations. 
 \end{simpassump}

\subsection{Id\`ele class characters}\label{Section:idele class characters}

\begin{defn} \label{def: Idele class character}
    An id\`ele class character
\[\chi=\sum_v \chi_v:\mathbb{A}_K^{\times}/K^{\times}\to \Q_p\]
is a \emph{continuous} homomorphism that decomposes as a sum of local characters $\chi_v$.     
\end{defn}

Below are some properties of id\`ele class characters:

\begin{itemize}
    \item For any prime $\mfq\nmid p$ we have $\chi_{\mfq}(O_{\mfq}^{\times})=0$ because of continuity. So if $\pi_{\mfq}$ is a uniformiser in $K_{\mfq}$,
then $\chi_{\mfq}$ is completely determined by $\chi_{\mfq}(\pi_{\mfq})$.
\item For any $\mfp\mid p$, there is a $\Q_p$
-linear map $t_{\mfp}^{\chi}$ such that we can decompose the local height as follows
\begin{equation} \label{diag: idele class char}
    \begin{tikzcd}
	{O_{\mfp}^{\times}} && {\Q_p} \\
	& {K_{\mfp}}
	\arrow["{\chi_{\mfp}}", from=1-1, to=1-3]
	\arrow["{\log_{\mfp}}"', from=1-1, to=2-2]
	\arrow["{t_{\mfp}^{\chi}}"', from=2-2, to=1-3]
\end{tikzcd}
\end{equation}
 because $\chi_{\mfp}$ take values in the torsion free group $(\Q_p,+)$ where $t_{\mfp}^{\chi}$ is a $\Q_p$-linear map.
\end{itemize}
If a continuous id\`ele class character $\chi$ is ramified at $\mfp$, that is, if the local character $\chi_{\mfp}$ does not
vanish on $O_{\mfp}^\times$, then we can extend $\log_{\mfp}:K_{\mfp}^\times\to K_{\mfp}$ in such a way that the diagram \eqref{diag: idele class char} remains commutative.

\begin{rmk}
    An id\`ele class character $\chi$ is admissible in the sense of \cite{mazur1983canonical} if the chosen abelian variety $A$ has good, ordinary reduction at places $v$ where $\chi_v$ is ramified, i.e. $\chi_v(O_v^{\times})\neq 0$. 
    If $\chi$ is admissible, then we can associate a pairing on $A(K)\times A(K)$.
\end{rmk}
\begin{example}
    If $K=\Q$, the unique id\`ele class character up to scalar multiplication is the  cyclotomic character $\chi\coloneqq \chi^{\text{cyc}}_{\Q}$. 
    If $x=(x_q)_q\in \A_{\Q}^{\times}$ then, $\chi^{\cyc}_{\Q}(x)=\log_p(x_p)-\sum_{q\neq p}\log_p(q^{v_q(x_q)})$, where $v_q$ denotes the $q$-adic valuation for $\Q_q^{\times}$. We choose $\log_p$ to be the Iwasawa branch of the $\log$ function, that is $\log_p(p)=0$. We note that $\chi(\Q^{\times})=0$, so $\chi$ factors through the id\`ele class group $\A_{\Q}^{\times}/\Q^{\times}$.
\end{example}

\begin{example} \label{ex:cyc character}
        For general $K$, let $\Nm_{K/\Q}:\A_K^{\times}\to \A_{\Q}^{\times}$ be induced by the usual norm function. The cyclotomic id\`ele class character is defined as $\chi_K^{\cyc}\coloneqq \chi_{\Q}^{\cyc}\circ\Nm_{K/\Q}$ 

    \end{example}

\begin{rmk}
    The character $\chi_K^{\cyc}$ is called `cyclotomic' because it corresponds to the Galois character which cuts out the cyclotomic $\Z_p$-extension via the identification of 
    $p$-adic Hecke characters with $p$-adic continuous Galois characters by class field theory.   
\end{rmk}
\begin{example}
    Suppose $K$ is an imaginary quadratic field. An anticyclotomic $p$-adic id\`ele class character is a continuous homomorphism
\[\chi : \A_K^{\times}/K^{\times}
\to\Z_p\]
such that $\chi\circ c = -\chi$. where $c$ is the map induced on $\A_K$ by complex conjugation.
\end{example}

\begin{rmk}
    By class field theory, the anticyclotomic character cuts out the anticyclotomic $\Z_p$-extension. 
\end{rmk}

\begin{lemma}\cite[Proposition 1.4] {balakrishnan2016shadow}\label{Lemma: Anticyc height}
Let $K$ be an imaginary quadratic field and let $K$ have class number $1$. Let $p$ be a prime which splits in $K$ and let $pO_K=\mfp_1\mfp_2$. Also let $\psi_1:K_{\mfp_1}\to \Q_p$ and $\psi_2:K_{\mfp_2}\to \Q_p$ be isomorphisms induced by the inclusion of $K$. Consider the map \begin{align*}
    \chi:\A_K^{\times}/K^{\times} &\to \Q_p\\
    (x_v)_v &\mapsto \log_p(\psi_1(\alpha)\cdot x_{\mfp_1})-\log_p(\psi_2(\alpha)\cdot x_{\mfp_2})
\end{align*}
where $\alpha\in K^{\times}$ is any element such that $\alpha x_v\in O_v^{\times}$
for all finite $v$. Then $\chi$ is the unique (up to scaling) non-trivial
anticyclotomic $p$-adic id\`ele class character.
\end{lemma}

\begin{proof}
    We refer the reader to \cite[Section 1]{balakrishnan2016shadow} for a proof, noting that the if $c$ is the involution on $A_K$ induced by complex conjugation, then $\psi_2(\alpha)=\psi_1(\alpha^c)$ for $\alpha\in K$. We note complex conjugation acts in the following manner (see for example \cite[Section V1.2]{neukirch2013algebraic})  
    \begin{align*}
        c:\A_K^{\times}\to  &\A_K^{\times}\\
            (x_v)_v\mapsto &(c(x_{v_c})).\\
    \end{align*}
     So, at the two places above $p$, complex conjugation acts by switching them. 
\end{proof}

 The restriction to class number one above is not essential, as there is an explicit formula (albeit more complicated) for class number not one. In this paper we look at number fields that have class number one, so will restrict our attention to this case.
\begin{lemma}\cite[Corollary 2.4]{balakrishnan2021explicit}
    Let $K$ have $2r_2$ complex places. Then the continuous id\`ele class characters of $K$ form a $\Q_p$-vector space
$V_K$ of dimension  $\geq r_2 + 1$. If Leopoldt’s conjecture \cite{Leopoldt1962} holds for $K$, then we
have $\dim_{\Q_p} V_K = r_2 + 1$. In particular, this holds if $K/\Q$ is an abelian extension.
\end{lemma}

For an imaginary quadratic field, the cyclotomic and anticyclotomic characters described above span the space of continuous id\`ele class characters.

\subsection{$p$-adic sigma function} \label{section: sigma function}
In their 1991 paper \cite{mazur1991p}, Mazur and Tate introduce the notion of $p$-adic sigma function, which we will recall in this section. 

Let $R$ be a complete discrete valuation ring with uniformiser $\pi$ and residue field
    $k=R/\pi R$ of characteristic $p>0$. Let $K$ be the field of fractions of $R$. Let $E/K$ be an elliptic curve over $K$. Let $E^f$ the formal group of $E$. We suppose that $E$ is ordinary, i.e. \[E^f_{\overline{k}}\cong \left({\G_m^f}\right)_{\overline{k}}.\] In that case there is a good $p$-adic analog of
the classical Weierstrass sigma function if $p$ is odd, and of its square if $p=2$. In \cite{mazur1991p}, the authors give 11 different characterisations of the $p$-adic sigma function, but we will present the one which is amenable to calculations. Let
\[x(t)=\frac{1}{t^2} +\cdots \in R((t))\]
be the formal Laurent series that expresses $x$ in terms of the local parameter
$t=-x/y$ at infinity. 

\begin{thm}[\protect{\cite[Theorem 1.3]{mazur2006computation}}\label{thm: p-adic sigma}]
    There is exactly one odd function $\sigma(t) = t+\cdots\in tR[[t]]$ and
constant $c\in R$ that together satisfy the differential equation 
\begin{equation}\label{eq:p-adic sigma}
x(t) + c = -\frac{d}{\omega} \left( \frac{1}{\sigma} \frac{d\sigma}{\omega}\right)    
\end{equation}
where $\omega$ is the invariant differential $dx/(2y + a_1x + a_3)$ associated with our
chosen Weierstrass equation for $E$.
\end{thm} 

An algorithm to compute the $p$-adic sigma function is given in \cite{mazur2006computation}. A more efficient algorithm to calculate it was given by David Harvey in \cite{harvey2008efficient}. This is the implementation we use. The implementation in SageMath \cite{sagemath} is for elliptic curves defined over the rationals. We modify the algorithm to use it for elliptic curves defined over (imaginary) quadratic extensions at split primes.

\begin{rmk} \label{rmk: quadraticity}

    If $m\in\Z$, and $Q\in E^f(\overline{R})$, then \cite[Theorem 3.1]{mazur1991p} states that\[\sigma(mQ)=\sigma(Q)^{m^2}f_m(Q).\]
where $f_m$ is the $m$th division polynomial relative to $\omega$. Now taking $\log$, we see \[\log(\sigma(mQ))=m^2\log(\sigma(Q))+\log(f_m(Q))\]

Informally, we say that, up to a correction term that the function $\log\sigma$ is quadratic.  
\end{rmk}


\subsection{Denominator} \label{sec: Denominator}
For any point $P=(x,y)\in E_v(K_v)$ which reduces to a non-singular point, for any finite place $v$, we can find $a_v,b_v,d_v$ which generate $O_v$, unique up to scaling by an element of $O_v^{\times}$  such that \[x=\frac{a_v}{d_v^2}, y=\frac{b_v}{d_v^3}\] since $O_v$ is a unique factorisation domain. 

Let $\pi_v$ be a uniformiser of $O_v$. Let $\val_v$ be the valuation on $K_v$. In \cite{mazur2006computation}, the authors define $\tilde{\sigma}_v(P)=\pi_v^{\val_v(d_v)}$ for $v\nmid p$. We note that $\tilde{\sigma}_v$ can be defined at places where the reduction \emph{is not} good and ordinary, whereas $\sigma_{v}$, the $v$-adic sigma function of $E_v$ cannot be defined at such places.

In \cite{mazur2006computation}, the authors say that $\tilde{\sigma}_{v}$ is a serviceable replacement for $\sigma_{v}$ in
the following sense. The $v$-adic valuation of $\tilde{\sigma}_v$ is the same as the $v$-adic valuation
of the $v$-adic sigma function. Since $\chi: \A_K^{\times}/K^{\times}\to \Q_p$ is unramified at $v\nmid p$ and consequently $\chi_v$ is 
sensitive only to the $v$-adic valuation of the input, this is an appropriate replacement for the sigma function. 

Note that for a global point $P=(x,y)\in E(K)$, we get points $P_v\in E(K_v)$ for all places $v$. In this case, Wuthrich \cite[Proposition 2]{wuthrich2004adic} shows that on a finite index subgroup $E^{\circ}(K)$ of  $E(K)$ we have a choice of a global denominator $d(P)\in O_K$ along with $a(P),b(P)\in O_K$ such that \[x=\frac{a(P)}{d(P)^2},y=\frac{b(P)}{d(P)^3},\] and $a(P),b(P),d(P)$ are pairwise coprime. Consequently $d_v=d(P)$ in $K_v$ up to $O_v^{\times}$ for all $v$. 

\begin{rmk} \label{rem: quadratic denominator}
    In fact in Wuthrich's thesis \cite[Proposition IV.3]{Wuthrichthesis}, he shows that for a point $P$ in $E^{\circ}(K)$, we have $d(mP)=d(P)^{m^2}f_m(P)$. Comparing this to Remark \ref{rmk: quadraticity}, we see the function \[h(P):=\log(\sigma(P))-\log(d(P))\] satisfies $h(mP)=m^2h(P)$ for $P\in E^{\circ}(K)$. Since the class group is trivial, $E^{\circ}(K)$ is precisely the subgroup of points that reduce to a non-singular point for all places $v$. 
\end{rmk}

\subsection{Computing heights for imaginary quadratic Fields}

In this section, we outline how one obtains heights from id\`ele class characters, denominators and $p$-adic sigma functions. We then give explicit formulas for the cyclotomic and anti-cyclotomic heights. 





For $\mfp\mid p$, let $\sigma_{\mfp}$ be the sigma function attached to the elliptic curve $E_{K_{\mfp}}$. Let $\mcE$ denote the minimal Weierstrass model. Note $\sigma_{\mfp}$ defines a function on the kernel of reduction $E_{\mfp}^1:=\ker (\mcE(K_{\mfp})\to \mcE(k_{\mfp})$) to $\Z_p$. Indeed, if $P=(x,y)\in E^1_{\mfp}$ and $t=-x/y$ then $\val_v(t)\geq 1$, and so $\sigma_{\mfp}(P)\coloneqq\sigma_{\mfp}(t)$ is a well-defined element of $O_{K_{\mfp}}$. If $m=\#\mcE(k_{\mfp})$, then any $P\in E(K),\, mP\in E^1_{\mfp}$.

For $v\nmid p$, if $P$ reduces to the identity component of the special fibre at $v$, $\tilde{\sigma}_v(P)=d_v$ as in Section \ref{sec: Denominator}. As discussed in Remark \ref{rem: quadratic denominator}, for points in $E^{\circ}(K)$, a finite index subgroup of $E(K)$, we have a global denominator, $d(P)\in O_K$. This is well defined up to an $O_K$ unit. If $n$ is the least common multiple of the Tamagawa numbers, then for any $P\in E(K)$, we have $nP\in E^{\circ}(K)$. 
\begin{defn}
    Let $E(K)_{\t{ht}}\subseteq E(K)$ be the set of non-torsion points $P$ which lie in $E^1_{\mfp}$ for all $\mfp\mid p$ and reduce to the identity component for $v\nmid p$. This is a finite index subgroup of $E(K)$.
\end{defn}

Fix $P\in E(K)_{\t{ht}}$ and an id\`ele class character $\chi$. For $\mfp\mid p$ and $v\nmid p$, we define \[h_{\mfp}^{\chi}(P)\coloneqq \chi_{\mfp}((\sigma_{\mfp}(P)))\hspace{1cm}\text{and} \hspace{1cm}h^{\chi}_{v}(P)\coloneqq \chi_{v}(\tilde{\sigma}_{v}(P)).\]

For $P\in E(K)$, let $n$ be such that $nP\in E(K)_{\t{ht}}$. As before, let $f_n$ be the $n$th division polynomial. For any place $v$, we extend local heights as follows:

\[h^{\chi}_{v}(P)=\frac{1}{n^2}\left(h^{\chi}_{v}(nP)-\chi_v(f_n(P))\right)\]

This definition does not depend on $n$; see Remark \ref{rmk: quadraticity} for $v\mid p$ and Remark \ref{rem: quadratic denominator} for $v\nmid p$. 

Finally, we have a formula for the global height 
\begin{equation}\label{Eqn: Ht formula}
    h^\chi(P)\coloneqq \sum_{v}h_v^{\chi}(P)=\frac{1}{n^2}\sum_v h_v^{\chi}(nP).
\end{equation}

 We now compute cyclotomic and anti-cyclotomic heights explicitly. We fix $K$, an imaginary quadratic field and $E/K$ an elliptic curve. Let $p$ be a split prime and write $\mfp_1,\mfp_2$ for the prime factors of $pO_K$. Let $\sigma_1,\sigma_2$ be the sigma functions attached to the curves $E_{K_{\mfp_1}}$ and $E_{K_{\mfp_2}}$ respectively. Fix a non-torsion point $P$, and choose an $n$ such that $nP\in E(K)_{\t{ht}}$. Let $nP=(x,y)$, and $t=-x/y$. Set $t_i=\psi_i(t)$ for $i=1,2$. Also, let $d$ be the denominator of $nP$ as in Section \ref{sec: Denominator}. 

\subsubsection{Cyclotomic Height}
Using \eqref{Eqn: Ht formula} and using the formula for the cyclotomic character in Example \ref{ex:cyc character}, we get
\[h^{\cyc}(P)=\frac{1}{n^2}\left(\log_p\left(\frac{\sigma_1(t_1)}{\sigma_2(t_2)}\right)+\log_p\left(\frac{\psi_1(d)}{\psi_2(d)}\right)\right).\]


\subsubsection{Anticyclotomic Height}

The formula to calculate with the anticyclotomic character is a bit more involved than the formula for the cyclotomic character. We prove a slight generalization of \cite[Proposition 2.4]{balakrishnan2016shadow} and use it to compute the anticyclotomic height. See Remark \ref{Rmk: Extension} for more details.

\begin{propn}[Anticyclotomic Height]\label{Propn: Anticyclotomic height}
Let $P \in E(K)_{\mathrm{ht}}$ be a non-torsion point. Then, the anticyclotomic $p$-adic
height of $P$, denoted $h^{\anti}(P)$ is given by the formula 

\[h^{\anti}(P)=\log_p\left(\frac{\sigma_1(t_1)}{\sigma_2(t_2)}\right)-\log_p\left(\frac{\psi_1(d)}{\psi_2(d)}\right)\]
    \end{propn}
\begin{proof}
Let $P=(x,y)$ reduce to the origin modulo $\mfp_1,\mfp_2$. By considering valuations, we see there exist $e_1,e_2\in \Z_{\geq 1}$ such that we have $\val_{\mfp_i}(x)=-2e_i,\val_{\mfp_i}(y)=-3e_i$ for $i=1,2$. We get \[\val_{\mfp_i}(t)=\val_{\mfp_i}(x)-\val_{\mfp_i}(y)=e_i.\] Furthermore $\val_{\mfp_i}(\sigma_i(t_i))=e_i$ since $\sigma_i(t)\in t+t^2\Z_p[t]$ for $i=1,2$. Let $\pi_1,\pi_2\in O_K$ generate $\mfp_1,\mfp_2$ respectively. Choosing $\alpha=\pi^{-e_1}$ in Lemma \ref{Lemma: Anticyc height}, we get\[\chi_{\mfp_1}^{\anti}\left(\sigma_1(P)\right)=\log_p(\psi_1(\pi_1^{-e_1}))+\log_p(\sigma_1(t_1))-\log_p\left(\psi_2(\pi_1^{-e_1})\right),\]
and choosing $\alpha=\pi_2^{-e_2}$, we get

\[\chi_{\mfp_2}^{\anti}\left(\sigma_2(P)\right)=\log_p(\psi_1(\pi_2^{-e_2}))-\log_p(\sigma_2(t_2))-\log_p\left(\psi_2(\pi_2^{-e_2})\right).\]
Therefore we get 
\begin{equation}\label{eqn: ht at p}
h_{\mfp_1}^{\anti}(P)+h^{\anti}_{\mfp_2}(P)=\log_p\left(\frac{\sigma_1(t_1)}{\sigma_2(t_2)}\right)-\log_p\left(\frac{\log_p(\psi_1(\pi_1^{e_1}\pi_2^{e_2}))}{\log_p(\psi_2(\pi_1^{e_1}\pi_2^{e_2}))}.\right)
\end{equation}

For local heights away from $p$ we get \[\sum_{v\nmid p\infty} h_v^{\anti}(P)=\sum_{v\nmid p\infty}\chi_v^{\anti}(d_v(P))=\sum_{v\nmid p\infty}\chi_v^{\anti}(d(P)) \]

by definition of $d(P)$. Choose $\beta \in O_K$ such that $d(P)=\pi_1^{e_1}\pi_2^{e_2}\beta$ and 
\[\val_{\mfp_1}(\beta)=\val_{\mfp_2}(\beta)=0.\]

Then by Lemma \ref{Lemma: Anticyc height} we get 

    \begin{align}
    \sum_{v\nmid p \infty}h^{\anti}_v(P)&=\log_p\left(\frac{\psi_1(\beta^{-1})}{\psi_2(\beta^{-1})}\right) \nonumber\\
    &= -\log_p\left(\frac{\psi_1(d(P))}{\psi_2(d(P))}\right)+\log_p\left(\frac{\psi_1(\pi_1^{e_1}\pi_2^{e_2})}{\psi_2(\pi_1^{e_1}\pi_2^{e_2})}\right) \label{eqn: anticyc away from p}.
\end{align}

Adding \eqref{eqn: ht at p} and \eqref{eqn: anticyc away from p}, we get the desired result.
\end{proof}

 \begin{rmk}\label{Rmk: Extension}
     Proposition \ref{Propn: Anticyclotomic height} is a mild generalization of Proposition 2.4 in \cite{balakrishnan2016shadow}. There are a few key differences:

     \begin{enumerate}
         \item Let $c$ denote the action of complex conjugation, which is the same as the non-trivial Galois automorphism of $K$.
         Since we are working with $E/K$, a curve which is not necessarily a base change of one over $\Q$, it is not necessarily the case that if $P=(x,y)\in E(K)$ then $(x^c,y^c)\in E(K)$. \item  We get two different sigma functions attached to $E$ above $p$.
     \end{enumerate}
      \end{rmk}
      
 We describe an algorithm to compute the cyclotomic and anti-cyclotomic heights.

\begin{algorithm}[Computing the Cyclotomic and Anticyclotomic Height] \label{algo:heights}
Input: An affine point $P=(x,y)\in E(K)$.\\
Output: The cyclotomic and anticyclotomic heights $h^{\cyc}(P)$ and $h^{\anti}(P)$.
\begin{enumerate}
    \item If $P$ is torsion, $h^{\cyc}(P)=h^{\anti}(P)=0$.
    \item We first compute the least common multiple of the Tamagawa numbers of the elliptic curve. Call this number $m$.
    \item Given $P\in E(K)$, we find the order of $P$ in $E_{\F_{\mfp_i}}$ for $i=1,2$. We call these $n_1,n_2$ and let $n=lcm(n_1,n_2,m)$. 
    \item Compute $nP$, and find a square root of the denominator of the $x$-coordinate $d(nP)$. Note we can do this since we have chosen $K$ to have class number $1$. 
    \item Let \[t=-\frac{x(nP)}{y(nP)},\] and let $t_1=\psi_1(t),t_2=\psi_2(t)$. Compute the sigma functions attached to $E_{K_{\mfp_i}}$ and call them $\sigma_i$ for $i=1,2$. Compute \[s_i\coloneqq \sigma_i(t_i)\] for $i=1,2$.
    \item Compute \[d_1=\psi_1(d(nP)),d_2=\psi_2(d(nP)).\]
    \item Return the heights:
    \[h^{\cyc}(P)=\frac{1}{n^2}\left(\log_p(s_1)+\log_p(s_2\right)-\frac{1}{n^2}\left(\log_p(d_1)+\log_p(d_2)\right)\]
    and 
     \[h^{\anti}(P)=\frac{1}{n^2}\left(\log_p(s_1)-\log_p(s_2\right)-\frac{1}{n^2}\left(\log_p(d_1)-\log_p(d_2)\right)\]
\end{enumerate}
\end{algorithm}

\section{Quadratic Chabauty for integral points of elliptic curves} \label{Sec:QC strategy}

Let $K$ be a quadratic imaginary number field with ring of integers $O_K$ and let $E/K$ be an elliptic curve. As before, let $K$ have class number one. Let $f(x,y)$ be the Weierstrass equation of the minimal Weierstrass model of $E$, and let \[\mcU\coloneqq \Spec O_K[x,y]/f(x,y). \]
We let $\mcE$ be the minimal regular resolution of projective closure of $\mcU$ in $\P^2_{O_K}$. We shall give a method to determine the set $\mcU(O_K)$, when $E(K)$ has rank $2$ and trivial torsion based on the proof of Theorem \ref{thm:Intergal points number fields}. 
If $E(K)$  has rank $0$, one just needs to compute the torsion points of $E$ to determine the set $\mcU(O_K)$. Since there are finitely many choices for the different groups $E(K)$ can be \cite{Kamienny:TorsionQuadratic}, this can be done relatively easily using division polynomials. 

In the case of $E(K)$ being rank one, one can generically find a $p$-adic linear functional and a height function which vanishes on the integral points of $\mcU$, as was done in \cite{balakrishnan2021explicit}. The case of rank $2$ is the one where one needs to genuinely use the two different height functions to determine the quadratic Chabauty set, and this is the case we shall explore. 
Let \[E(K)=\Z\cdot P\oplus \Z \cdot Q.\] 
Fix an id\`ele class character $\chi$.
In Theorem \ref{thm:Intergal points number fields}, we have a function $\rho^{\chi}$ for a choice of id\`ele class character $\chi$, and a set $T^{\chi}$ such that $\rho^{\chi}(\psi(\mcU(O_K)))\subseteq T^{\chi}$. We shall explain how to obtain the function $\rho^{\chi}$ and the set $T^{\chi}$ using the $p$-adic heights considered in Section \ref{sec: Heights}.

We first recall some preliminaries on Coleman integration.

\subsection{Coleman integration and height functions}
Let $L/\Q_p$ be a finite extension and let $X_{an}/L$ be a rigid analytic curve. See \cite[Section 1]{coleman1985torsion} for notions on rigid geometry used in the following theorem. See \cite{fresnel2012rigid} for a more thorough introduction to rigid geometry.

\begin{thm}[Coleman, \protect{\cite[Proposition 2.4, Theorem 2.7]{coleman1985torsion}}]\label{Thm:Integral-properties}
    Let $\eta,\xi$ be $1$-forms on a wide open V of $X_{\text{an}}$ and
$A, B, C \in V (L)$. Let $a, b \in L.$ The definite Coleman integral has the following properties:

\begin{enumerate}
    \item \textbf{Linearity}: \[\int_A^B(a\eta + b\xi) = a \int_A^B \eta + b \int_A^B \xi.\]
    \item \textbf{Additivity in endpoints}:\[ \int_A^B \xi = \int^B_C \xi + \int^C_A \xi\]
    \item \textbf{Change of variables}: If 
$\phi : X_{\mathrm{an}} \to X_{\mathrm{an}}'$ is a rigid analytic map then \[\int^B_A \phi^*\xi = \int^{\phi(B)}_{\phi(A)} \xi\]
\item \textbf{Fundamental theorem of calculus}: \[\int_A^B df = f(B) - f(A)\] for $f$ a rigid analytic function on
$X$.
\end{enumerate}
\end{thm}
We now state an easy consequence of the properties of the Coleman integral, which we will use.
\begin{lemma}[\protect{\cite[Theorem 2.8]{coleman1985torsion}}] \label{lemma:addition on EC points}

Let $E/L$ be an elliptic curve and $P_1,P_2\in E(L)$ be points on the elliptic curve. 
Let $\omega\in H^0(E,\Omega^1)$ be an invariant differential. Then 
\[\int^{P_1+P_2}_{\mcO}\omega=\int^{P_1}_{\mcO}\omega+\int^{P_2}_{\mcO}\omega\]
\end{lemma}

We use the implementation of Coleman integration in \cite{balakrishnan2010explicit}.

    


 Let $p$ be a split prime and write $\mfp_1,\mfp_2$ for the prime factors of $pO_K $. Further, we let $\psi_1,\psi_2:K\into \Q_p$ be the natural embeddings. Let \begin{align*}
    &E_{K_{\mfp_1}}\coloneqq E\times_{\psi_1}\Spec\Q_p, \hspace{0.5 cm} E_{K_{\mfp_2}}\coloneqq E\times_{\psi_2}\Spec\Q_p.
\end{align*}

Fix $\omega\in H^0(E,\Omega^1)$, a choice of invariant differential, and let $\omega_i=\psi_i^*\omega$ be the pullbacks to $E_{K_{\mfp_i}}$ for $i=1,2$. For $P_1\in E_{K_{\mfp_1}}(\Q_p),P_2 \in E_{\K_{\mfp_2}}(\Q_p)$ we let 
\[f_1(P_1)\coloneqq \int^{P_1}_{O}\omega_1 \hspace{ 2cm}\text{and} \hspace{2 cm} f_2(P_2)\coloneqq \int^{P_2}_O \omega_2\]
where the integrals are Coleman integrals. In particular for $R\in E(K)$ we let 
\[f_1(R)\coloneqq \int^{\psi_1(R)}_O \omega_1 \hspace{2 cm}\text{and}\hspace{2cm} f_2(R)\coloneqq\int^{\psi_2(R)}_O \omega_2.\]
By Lemma \ref{lemma:addition on EC points}, $f_i$ are linear functions on $E_{K_{\mfp_i}}(\Q_p)$ for $i=1,2$. Consider \[\t{Span}_{\Q_p}\{f_1,f_2\} \subseteq (E(K)\otimes \Q_p)^{\vee}\]
where $(E(K)\otimes \Q_p)^{\vee}=\Hom_{\Q_p}(E(K)\otimes\Q_p,\Q_p)$. 

If $f_1,f_2$ are linearly independent, this inclusion is an equality (since $E(K)$ has rank $2$), and we can use the idea of quadratic Chabauty and restriction of scalars to find the integral points $\mcU(O_K)$. We now assume the following condition holds:

\begin{cond} \label{cond: lin independence}
    $\t{Span}_{\Q_p}\{f_1,f_2\}=(E(K)\otimes \Q_p)^{\vee}$.
\end{cond}
If Condition \ref{cond: lin independence}, one can compute a linear combination of $f_1$ and $f_2$ which vanishes on the integral points of $E$ via a version of Siksek's method. 

\begin{defn}\label{Defn: Bilinear forms}
For $R_1\in E_{K_{\mfp_1}}(\Q_p),R_2\in E_{K_{\mfp_2}}(\Q_p)$, let \[g_{ij}(R_1,R_2)=\frac{1}{2}(f_i(R_1)f_j(R_2)+f_i(R_1)f_j(R_2))\] 
for $1\leq i\leq j\leq 2$. 
\end{defn}
\begin{rmk}
When Condition \ref{cond: lin independence} holds $\t{Span}_{\Q_p}\{g_{11},g_{12},g_{22}\}$ is the full space of symmetric $\Q_p$-bilinear forms on $E(K)\otimes \Q_p$. It is essential in computing locally analytic expansions of the functions $\rho^{\chi}$  
\end{rmk}
\subsection{Finding the quadratic Chabauty sets for elliptic curves} \label{Subsection: QC for elliptic curves}
The global height pairing is bilinear, so if Condition \ref{cond: lin independence} holds, it can be expressed as a linear combination
\begin{equation}\label{eqn: hts and coleman integral}
h^{\chi}=\alpha^{\chi}_{11}g_{11}+\alpha^{\chi}_{12}g_{12}+\alpha^{\chi}_{22}g_{22},
\end{equation}
where the coefficients $\alpha_{i,j}\in Q_p$.
Now we return to sketch an algorithm to find a set $B_p\subseteq E(K\otimes\Q_p)$ which contains the image $\psi_1(\mcU(O_K))\times \psi_2(\mcU(O_K))$. We need the following properties of heights:

\begin{enumerate}
    \item For $v\nmid p$, $h_v^{\chi}$ takes only finitely many values on $O_K$-points. Call this set of values $T_v^{\chi}$. If $v$ does not divide the discriminant of $E$, then $T_v^{\chi}=\{0\}$.
    This follows from the fact that the local heights of Coleman--Gross \cite{coleman1989p} coincide with the height of Mazur and Tate by work of Balakrishnan and Besser \cite[Corollary 4.3]{balakrishnan2015coleman}. The local heights of Coleman and Gross away from places above $p$ are given by
    \[h_v(R)=(R,R)\chi(v)\]
    where $(\cdot,\cdot)$ is an intersection pairing. In \cite[Lemma 2.4]{balakrishnan2017computing}, the authors show the intersection pairing, and hence the local height at $v\nmid p$ only depends on the connected component that $R$ reduces onto in the special fibre. Hence, this set is finite and in fact equal to $\{0\}$ for almost all $v$.

    To compute this set for elliptic curves, we follow the algorithm in \cite[Section 5]{silverman1988computing} and use
    \cite[Proposition 6]{cremona2006height}.
    
\item  The global height $h^{\chi}$ and local height $h^{\chi}_{\mfp}$ for $\mfp\mid p$ can be expressed as locally analytic $p$-adic functions in each residue disk and computed to desired $p$-adic and $t_1,t_2$-adic precision as elements of $\Z_p[[t_1,t_2]]$. 
For a point $P$ such that $mP$ lies in an appropriate neighbourhood of the identity, we use the identity $h^{\chi}_{\mfp}(P)= \frac{1}{m^2}\left(h_{\mfp}^{\chi}(mP)-\chi_{|mfp}(f_m(P))\right)$

To compute analytic expansions of $h^{\chi}$ in all residue disks we first find $\alpha_{11},\alpha_{12},\alpha_{22}$ in Equation \ref{eqn: hts and coleman integral} by plugging in $P,Q,P+Q$. Then, we can compute Coleman integrals by finding a local parametre at a fixed point in each disk and using it to compute a tiny integral.


Hence we can compute 

\[\rho^{\chi}\coloneqq h^{\chi}-\sum_{\mfp|p}h_{\mfp}^{\chi}\]

as a $p$-adic analytic function up to desired $p$-adic precision in each residue disk.
\end{enumerate}

Thus $\rho^{\chi}$ is a $p$-adic analytic function in two variables, and \[\rho^{\chi}(\psi_1(\mcU(O_K))\times \psi_2(\mcU(O_K)))\subseteq T^{\chi},\] 
where $T^{\chi}$ is determined by the sets $T_v^{\chi}$ above.

For $K$ imaginary quadratic, as discussed in Section \ref{Section:idele class characters}, one has two linearly independent characters $\chi^{\cyc}$, the cyclotomic character and $\chi^{\anti}$ the anticyclotomic character. Let $\bar{\rho}=(\rho^{\cyc},\rho^{\anti})$. We wish to solve the system of equations $\bar{\rho}=(t_1,t_2)$ for $(t_1,t_2)\in T^{\cyc}\times T^{\anti}$.

This can be done using a multivariate Hensel's lemma. We can lift mod $p^n$ solutions using a Newton's method type argument uniquely in disks of small radius under certain conditions satisfied by the Jacobian matrix of partials of $\bar{\rho}$. See, for example \cite[Appendix A]{balakrishnan2021explicit}. Since finitely many disks cover $\mcU(\Z_{\mfp_i})$, one gets finitely many $\mcU(\Z_{\mfp_1})\times\mcU(\Z_{\mfp_2})$-solutions $(z_1,z_2)$ such that $\bar{\rho}(z_1,z_2)=(t_1,t_2)$ for $(t_1,t_2)\in T^{\cyc}\times T^{\anti}$.
Some adjustments need to be made when the Jacobian matrix does not satisfy the required criterion, which is also discussed in \cite[Appendix A]{balakrishnan2021explicit}.

We call the set of solutions obtained above the \emph{quadratic Chabauty set} and denote it $B_p\subseteq \mcU(O_K\otimes \Z_p)$. We use a version of algdep which relies on the LLL \cite{lenstra1982factoring} algorithm to recognise $O_K$-points for $p$-adic points. We call the set of points which we have not been able to identify as \emph{mock integral points} and use the notation $A_p$ to identify them. We would like to use some sieving method to show that these points are not integral, but unfortunately the Mordell--Weil sieve can not be used for elliptic curves. 

We discuss an alternative sieve in Section \ref{sec: sieve}. 
Here is an overview of the algorithm to be followed to find the set $A_p$ above.
\begin{algorithm}[Finding the quadratic Chabauty set for an elliptic curve]
\label{algo: QC}\ \\
Input: An elliptic curve $E/K$ over an imaginary quadratic field such that $E(K)=Z\cdot P\oplus \Z\cdot\Q$. A prime $p$ which splits in $K$ and $E$ has good, ordinary reduction at both primes above $p$.

Output: A finite set $B_p\subseteq \mcU(O_K\otimes \Z_p)$ containing $\sigma(\mcU(O_K))$ and a set $A_p$ containing mock integral points.
\begin{enumerate}
    \item Find the Mordell--Weil generators of $E(K)$ and let $E(K)=\Z\cdot P\oplus \Z\cdot Q$. To speed up the algorithm, one can optionally fix a positive integer $C$ and compute
\[E(O_K)_{\text{known}} = \left\{nP + mQ : |n|, |m| \leq C \text{ and }nQ + mR \in \mcU(O_K)\right\}.\]
\item Calculate the sets $T^{\cyc}$ and $T^{\anti}$ by computing the possible values of intersection pairings as in \cite{cremona2006height} and \cite{silverman1988computing}.
\item For $\chi=\chi^{\cyc}, \chi=\chi^{\anti}$ find an analytic expansion of $h^{\chi}-\sum_{\mfp|p}h^{\chi}_v$ in all residue disks. This will be done in two steps.
\begin{enumerate}
    \item Compute the local height at $\mfp \mid p$ by using the $p$-adic sigma function as in Section \ref{section: sigma function}. 
    \item Consider the system \[\begin{pmatrix}h(P,P)\\
    h(P,Q)\\
    h(Q,Q)\end{pmatrix}=
    \begin{pmatrix}
        g_{11}(P,P) & g_{12}(P,P) & g_{22}(P,P)\\
        g_{11}(P,Q) & g_{12}(P,Q) & g_{22}(P,Q)\\
        g_{11}(Q,Q) & g_{12}(Q,Q) & g_{22}(Q,Q)
    \end{pmatrix}
    \begin{pmatrix}
        \alpha_{11}\\ \alpha_{12}\\ \alpha_{22}
    \end{pmatrix}\]
    and solve for $\alpha_{11},\alpha_{12},\alpha_{22}$. Using implementations of Coleman integration for the matrix on the right and Algorithm \ref{algo:heights} for the heights on the left, we can find $\alpha_{11},\alpha_{12},\alpha_{22}$. Hence we can express $h-\sum_{\mfp|p}h_v$ as a $p$-adic power series in two variables.
\end{enumerate}
    \item Let $\rho^{\chi}=h^{\chi}-\sum_{\mfp|p}h^{\chi}_{\mfp}$.We solve the system $(\rho^{\chi^{\cyc}},\rho^{\chi^{\anti}})=(t_1,t_2)$ for all $(t_1,t_2)\in T^{\cyc}\times T^{\anti}$ using multivariate Hensel's lemma. Pull back solutions to points in $\mcU(O_K)$. This is the set $B_p$. 
    
    \item We identify elements of $B_p$ which lie in $\mcU(O_K)$ using $p$-adic LLL. We let $A_p$ be the complement of the points which are recognised as integral. 

\end{enumerate}

\end{algorithm}

\begin{rmk}
    The bottleneck in this algorithm is solving the system of $p$-adic power series. One would like to optimise solving these multivariate series to make this algorithm faster. 
\end{rmk}

In the next section, we discuss a method that shows certain points in $A_p$ are not $O_K$-points.



\section{A sieve for integral points on elliptic curves} \label{sec: sieve}

As before, let $p$ be a prime which splits in the chosen imaginary quadratic field $K$. On running the quadratic Chabauty algorithm for an elliptic curve to find integral points, one gets the following output:

\[B_p\coloneqq\left\{(R_1,R_2)\in \mathcal{U}(\Z_{p_1})\times\mathcal{U}(\Z_{p_2}):\rho^{\cyc}(R_1,R_2)\in T^{\cyc},\rho^{\anti}(R_1,R_2)\in T^{\anti}\right\}.\]

We often have $p$-adic points which are not recognised as algebraic points. For a curve with genus $g\geq 2$, one uses the Mordell--Weil sieve to show auxiliary $p$-adic points are not algebraic, but this is not applicable in our case.  

Instead, we outline a method that first appeared in the appendix of the paper \cite{balakrishnan2017computing}.
We first choose two primes $p,q$ for which we run the quadratic Chabauty algorithm. As always, we assume the primes split in the chosen imaginary quadratic field,  say $p O_K=\mfp _1\mfp _2$ and $q O_K=\mfq_1\mfq_2$. Also assume $E$ has good reduction at $\mfp_1,\mfp_2,\mfq_1,\mfq_2$. The primes must satisfy the following condition:

\begin{cond}\label{condition: reduction}
    Let $\F_{\mfp_i}\coloneqq O_K/\mfp_i$, $\F_{\mfq_i}\coloneqq O_K/\mfq_i$. We require $q\mid \#\mcE(\F_{\mfp_i})$ for at least one of $i=1,2$, \emph{and} $p\mid \#\mcE(\F_{\mfq_j})$ for at least one of $j=1,2$. This is equivalent to $p\mid \#\mcE(O_K/q)$ \emph{and} $q\mid \#\mcE(O_K/p)$ since $\mcE(O_K/p)=\mcE(\F_{\mfp_1})\times \mcE(\F_{\mfp_2})$ and $\mcE(O_K/q)=\mcE(\F_{\mfq_1})\times \mcE(\F_{\mfq_2})$.
\end{cond} 
Since we have chosen a curve such that $$E(K)=\Z\cdot P \oplus \Z\cdot Q,$$
we know any integral point, $R$ can be expressed as \begin{equation}\label{eqn: OK-point}
    R=mP+nQ
\end{equation}
for some unique integers $m,n$.

We shall find constraints on the $m$ and $n$ mod $p$ and mod $q$ using Coleman integrals (log information) and using the cardinalities of the reductions at $p,q$ which satisfy Condition \ref{condition: reduction}. We discuss this in the upcoming sections. For an abelian group $G$ and a positive integer $n$, we let $G/n$ denote the quotient $G/nG$. 
\subsection{Reduction information}\label{sec:reduction}

Let $p,q$ be rational primes which satisfy Condition \ref{condition: reduction}. Consider the following diagram:
\begin{figure}[h]
    \centering
\begin{tikzcd}
\mcU(O_K) \arrow[r] \arrow[d] & B_p \arrow[d]             \\
\mcE(O_K) \arrow[r] \arrow[d, , ]     & E(K\otimes\Q_p) \arrow[d ] \\
                           \mcE(O_K)/q \arrow[r ,"\red"]   & \mcE(O_K/p)/q  .          
\end{tikzcd}
\caption{Reduction data.}
\label{fig:reddata}
    \end{figure}

Let $\bar{\mcR}$ denote the image of $\mcR$ under the map $B_p\to E(K\otimes\Q_p)\to \mcE(O_K/p)/q$.
\begin{defn}\label{D:modqredn}
Let $\red_{\mcR}\coloneqq \red^{-1}(\bar{\mcR}) \subset \mcE(O_K)/q$. 
We call this set the collection of mod $q$ reduction constraints on a putative integral point corresponding to the $p$-adic point $\mcR$. 
\end{defn}


Since we have fixed generators $P,Q$, we have an isomorphism $\mcE(O_K)\cong \Z^2$, and with this choice of basis, we obtain (mod $q$) constraints on the coordinates $m$ and $n$ of a putative integral point giving rise to $\mcR$ in this basis. 

\begin{rmk}
By condition \ref{condition: reduction}, the $\F_p$-vector space $\mcE(O_K/p)/q$ is not isomorphic to \{0\}, so we do get meaningful information by considering $\red_{\mcR}$.
\end{rmk}  

Let $\bar{P}_i,\bar{Q}_i$ for $i=1,2$ denote the image of the generators in $E_{\F_{\mfq_i}}(\F_q)$ for $i=1,2$, and let $\bar{R}_i$ denote the reduction of $\mcR=(R_1,R_2)\in B_p$. Then we can compute $\red_{\mcR}$ as

\[ \mathrm{red}_{\mcR}=\left\{(m,n)\mid m,n\in \Z/p, \bar{R}_1=m\bar{P}_1+n\bar{Q}_1,\bar{R}_2=m\bar{P}_2+n\bar{Q}_2 
\right\}.\]

If this set is empty, we can immediately discard the point $(R_1,R_2)$ as not coming from a point $R\in \mcE(O_K)$, by the commutativity of Figure \ref{fig:reddata}. Else we record $\red_{\mcR}$.


\subsection{Log information}
\label{section: log}

Recall that we have a map $\psi: E(K)\to E(K\otimes\Q_p),P\mapsto (\psi(P))$ which embeds points diagonally. Let $T_1,T_2$ denote the tangent spaces of $E(K_{\mfp_1})$ and $E(K_{\mfp_2})$. For $i=1,2$, we had $f_1,f_2$ in Condition \ref{cond: lin independence} as functionals on $H^0(E_{\Q_{\mfp_i}},\Omega^1)\cong T_i^*$. We have a natural log map \begin{align*}
    \log:E(K\otimes\Q_p) &\to T_1\times T_2\\
    (R_1,R_2) &\mapsto (f_1(R_1),f_2(R_2))
\end{align*}


  
We have the following commutative diagram which encapsulates the objects and maps we have discussed:

\begin{center}
\begin{figure}[h]
    \centering

\begin{tikzcd}
\mcU(O_K) \arrow[r] \arrow[d]                                           & B_p \arrow[r] & E(K\otimes\Q_p) \arrow[d, "\log"] \\
\mcE(K) \arrow[rr, "\log|_{\mcE(K)}"] \arrow[d, "\mathrm{red}"] &               & T_1\times T_2                     \\
\mcE(K)/p                                                               &               &                                  
\end{tikzcd}\caption{Log data}
\label{fig:logdata}
    \end{figure}

\end{center}


\begin{defn}\label{D:loginf}
Let $\log_{\mcR}\coloneqq \t{red}
(\log|_{\mcE(K)}^{-1}{\mcR}))\subseteq \mcE(K)/p$.  We call this set the collection of mod $p$ log constraints on a putative integral point corresponding to the $p$-adic point $\mcR$.
\end{defn}
Suppose $\mcR=(R_1,R_2)\in A_p$. To calculate $\log_{\mcR}$, we solve for $m,n$ in the following equations: 

    \begin{align}\label{Eqn: log equation}
    \begin{split}
f_1(R_1)  &=  m f_1(P)+n f_1(Q)\\
f_2(R_2) &=  m f_2(P)+n f_2(Q).
\end{split}
\end{align}
If $M$ is invertible, we let \[M=\begin{pmatrix}
 f_1(P_1) & f_1(Q_1)\\
 f_2(P_2) & f_2(Q_2)
\end{pmatrix},\]
and let \begin{equation}\label{eqn:log info}
    \begin{pmatrix}
    m\\n
\end{pmatrix}=M^{-1}\begin{pmatrix}
    f_1(R_1)\\f_2(R_2)
\end{pmatrix},
\end{equation}
where $m,n\in\Q_p$. If $m,n\not\in\Z_p$, then the point $(R_1,R_2)$ is not integral. If $m,n\in \Z_p$, we can take mod $p$ reductions to get $\log_{\mcR}$. If $M$ is not invertible, then we need to check if \eqref{eqn:log info} has infinitely many solutions or no solutions. If it has no solution, we can discard the point $\mcR$, else, we reduce mod $p$ and record $\log_{\mcR}$.





\subsection{Comparing log and reduction information}

For $\mcR \in B_p$, let $\red_{\mcR} \times \log_{\mcR} \subset \mc{E}(O_K)/q \times (O_K)/p$ be as in Definition~\ref{D:modqredn} and Definition~\ref{D:loginf}. Similarly, for $\mcS \in B_q$, let $\log_{\mcS} \times \red_{\mcS} \subset \mcE(O_K)/q \times \mcE(O_K)/p$. Let $\mathfrak{S}$ be the union of  $\log_{\mcS} \times \red_{\mcS}$ for all $\mcS$ in $B_q$. 
\begin{lemma}\label{L:discard}
 Assume condition~\ref{condition: reduction}. Let $\mcR \in B_p$.   If $(\red_R \times \log_R) \cap \mathfrak{S}$ is empty, then $\mcR$ is not the image of a point in $\mcU(O_K)$.
\end{lemma}
\begin{proof}
By the commutativity of Diagrams~\ref{fig:reddata} and \ref{fig:logdata} , it follows that the image of a global point in $\mcU(O_K)$ in $\mcE(K)/q \times \mcE(K)/p$ must be in the intersection of $(\red_R \times \log_R) \cap \mathfrak{S}$.
\end{proof}

We now list the steps to carry out the sieve. 
\begin{algorithm}\label{algo:sieve }
Input: Sets $A_p\subseteq \mcE(O_K\otimes\Z_p)$ and $A_q\subseteq \mcE(O_K\otimes\Z_q)$.\\
Output: Subsets $A'_p\subset A_p$  and $A'_q\subset A_q$, which only contain points which satisfy reduction and log considerations.
\begin{enumerate}
\item For all $\mcR\in A_p$, compute $\red_{\mcR}$. Let $R_q=\{\red_\mcR:\mcR\in A_p\}$. Let $R_p$ be the analogous set for the set $A_q$. 
\item For $\mcR\in A_p$, calculate $\log_{\mcR}$ as in \ref{section: log}. 
Let $L_p\coloneqq \{\log_{\mcR}:\mcR\in A_p\}$. Similarly, let $L_q$ denote the analogous set for $A_q$.
\item For each $\mcR$, compute 
\begin{equation}\label{Eqn: logred}
    \left(\log_{\mcR}\times \red_{\mcR}\right)\bigcap\left(\bigcup_{\mcS\in A_q}\red_{\mcS}\times\log_{\mcS}\right)
\end{equation}
If the set in \eqref{Eqn: logred} is empty, then we discard $\mcR$. 
\item If the set in \eqref{Eqn: logred} is not empty, we record each $\mcS\in A_q$ such that 
$\left(\log_{\mcR}\times \red_{\mcR}\right)\cap \left(\log_{\mcS}\times \red_{\mcS}\right)\neq \emptyset$. We check for each such $\mcS$ if the sum of local heights away from $q$ corresponds to the sum of local height away from $p$ computed for $\mcR$. If they correspond to different local heights, we discard $\mcR$. Else, we append $\mcR$ to $A'_p$.
\item We run through Step 1 to Step 4 for all $\mcR\in A_p$. This leaves us with a subset $A'_p\subseteq A_p$ of points that satisfy the restrictions obtained from log information and reduction information. We also get a similar set $A'_q\subseteq A_q$.
\end{enumerate}
\end{algorithm}

\section{Examples} \label{Sec:Examples}

Let $K=\Q(\zeta_3)$, and $O_K=\Z[\zeta_3]$.
    
\begin{example}\label{ex: First computation}
    Consider the scheme $\mcU_1/O_K$ cut out by the Weierstrass equation 
     \[y^2+(\zeta_3+2)y=x^3+(-\zeta_3-2)x^2+(\zeta_3+1)x.\]

We let $E$ be the projectivisation of the generic fibre of $\mcU_1$.
Let $pO_K=\mfp_{1}\mfp_2$ for $p=7$ and $qO_K=\mfq_1\mfq_2$ for $q=13$.
We have the following information about this curve:

\begin{enumerate}
    \item  $E(K)= \Z\cdot P \oplus \Z\cdot Q$ where $P=(1,0),Q=(\zeta_3+1,0)$.
    \item $E_{\F_{\mfp_1}}(\F_p)\cong \Z/13\Z\cong E_{\F_{\mfp_2}}(\F_p)$.
    \item $E_{\F_{\mfq_1}}(\F_q)\cong \Z/7\Z, E_{\F_{\mfq_2}} (\F_q)\cong\Z/19\Z$.
    \item The curve has bad reduction at the single prime $\mfp=(-22\zeta_3-9)$ and $\Nm(\mfp)=367$. It has Kodaira Symbol $II$, Tamagawa number $1$.  The discriminant and conductor of the curve are both $\mfp^2$. For more information about this curve, see the LMFDB page at \url{https://www.lmfdb.org/EllipticCurve/2.0.3.1/134689.3/CMa/1}.
    \item The sets $T^{\cyc}$ and $T^{\anti}$ as described in Algorithm \ref{algo: QC} are $\{0\}$ for $i=1,2$.
\end{enumerate}
We first search for points of small height and get a list of $12$ $O_K$-points. We wish to certify these are the only points of $\mcU_1$.

On finding the quadratic Chabauty sets for $p=7,q=13$ since they satisfy Condition \ref{condition: reduction}, we get $\#A_p=204,\#A_q=108$. On applying the sieve described in Algorithm \ref{algo:sieve }, all the classes are eliminated, that is all the extra $7$-adic points and $13$-adic points are eliminated via mod $7$ and mod $13$ considerations. Hence we get 
\begin{align*}
    \mcU_1(O_K)=&\{(-3 : -8\zeta_3 - 4), (-3 : 7\zeta_3 +2 ),(4\zeta_3+4 : -8\zeta_3 - 4), (0 : 0 ), (0 : -\zeta_3 - 2), (\zeta_3+1 : 0 ),\\ (\zeta_3+1 : -\zeta_3 - 2),
    & (4\zeta_3+4 : 7\zeta_3 +2 ), (1 : 0), (1 : -\zeta_3 - 2), (-3\zeta_3 + 1 : -8\zeta_3 - 4), (-3\zeta_3 + 1 : 7\zeta_3 + 2)\}.
    \end{align*} 
\end{example}

\begin{example}
    Consider the scheme $\mcU_2/O_K$ cut out by the Weierstrass equation
    \[y^2+(-\zeta_3+1)y=x^3+(\zeta_3-1)x^2-\zeta_3x\]

    More infromation about this curve can be found at \href{}{}This is the Galois-conjugate of the scheme $\mcU_1$ in Example \ref{ex: First computation}.  We have a bijection of $O_K$-points via the non-trivial Galois automorphism of $\Q(\zeta_3)$; that is, there is a bijection 

    \[\mcU_1(O_K)\longleftrightarrow
   \mcU_2(O_K).\]

   Hence we get $\#\mcU_2(O_K)=12$.  
   
\end{example}

\begin{example}
    Consider the scheme $\mcU_3\subseteq \A^2_{O_K}$ given by 
    \[y^2 + (\zeta_3+1)y = x^3 + (-\zeta_3-2)x^2 + (\zeta_3+1)x + (\zeta_3+2)\]

There are at least $24$ points of small height in $\mcU_3$, and we compute a finite superset $\mcU_3(\Z_p)_{\mathrm{sieve}}\subseteq \mcU(O_K\otimes\Z_p)$ of the integral points. 
We let $E$ be the elliptic curve given by the Weierstrass equation above. For more information about this curve, see the LMFDB page at \\
\url{https://www.lmfdb.org/EllipticCurve/2.0.3.1/47089.9/CMa/1}.

Let $pO_K=\mfp_{1}\mfp_2$ for $p=13$ and $qO_K=\mfq_1\mfq_2$ for $q=19$.
We have the following information about this curve:

\begin{itemize}
    \item  $E(K)= \Z\cdot P\oplus \Z\cdot Q$ where $P=(-1,-2\zeta_3-1),Q=(2\zeta_3+2,\zeta_3)$.
    \item $E_{\F_{\mfp_1}}(\F_p)\cong \Z/19\Z, E_{\F_{\mfp_2}}(\F_p)\cong \Z/21\Z$.
    \item $E_{\F_{\mfq_1}}(\F_q)\cong \Z/13\Z\cong E_{\F_{\mfq_2}} (\F_q)$.
    \item The curve has bad reduction at a prime above $7$, $v_1=(3a-2)$ and at a prime above $31$, $v_2=(6a-5)$. Over $v_1$, the Kodaira Symbol is $IV$, with Tamagawa number $3$ and at $v_2$, the Kodaira symbol is $II$ with Tamagawa number $2$.  The conductor is $v_1^2v_2^2$ and the discriminant is $v_1^4v_2^2$.
    \item The sets $T^{\cyc}=\left\{0,-\frac{2}{3}\chi^{\cyc}(v_2)\right\}$ and $T^{\anti}=\left\{0,-\frac{2}{3}\chi^{\anti}(v_2)\right\}$ as described in Algorithm \ref{algo: QC}. 
    The set $A_p$ has cardinality $672$ while the set $A_q$ has cardinality $216$. Performing the sieve described in Algorithm \ref{algo:sieve }, we have $24$ points which are not accounted for. So we have a set $\mcU_3(\Z_p)_{\mathrm{sieve}}$ such that \[\mcU_3(O_K)_{\mathrm{known}}\subseteq\mcU_3(O_K)\subseteq \mcU_3(\Z_p)_{\mathrm{sieve}}\subseteq \mcU(O_K\otimes \Z_p)\]
     such that $\#\mcU_3(O_K)_{\mathrm{known}}=24$ and $\mcU_3(\Z_p)_{\mathrm{sieve}}=48$.
\end{itemize}

\end{example}


%
\bibliographystyle{amsalpha}

\bibliography{bibliography}

\end{document}